\documentclass[12pt,a4]{article} 
\usepackage{amsmath,amssymb,eucal} 
\usepackage{bbm}  
\font\tenrm=cmr10

\font\cmssl=cmss10 at 12 pt  
\font\bigss=cmssdc10 scaled 2300

\font\cmsslll=cmss10 at 14 pt

\renewcommand{\a}{\alpha}

\renewcommand{\o}{\omega}

\newcommand{\bC}{\mathbb{C}}  
\newcommand{\bR}{\mathbb{R}}

\newcommand{\id}   {{\mathbbm{1}}}

\renewcommand{\square}{\kern1pt\vbox  
               {\hrule height 0.6pt\hbox{\vrule width 0.6pt\hskip 3pt  
    \vbox{\vskip 6pt}\hskip 3pt\vrule width 0.6pt}\hrule height0.6pt}  
    \kern1pt}

\newtheorem{Th}{Theorem}  
\newtheorem{Prop}{Proposition}  
\newtheorem{Cor}{Corollary}  
\newtheorem{Lem}{Lemma}  
\newtheorem{Def}{Definition}  
  
\newcommand{\bt}{\begin{Th}\ \ }  
\newcommand{\et}{\end{Th}}  
\newcommand{\bp}{\begin{Prop}\ \ }  
\newcommand{\ep}{\end{Prop}}  
\newcommand{\bc}{\begin{Cor}\ \ }  
\newcommand{\ec}{\end{Cor}}  
\newcommand{\bl}{\begin{Lem}\ \ }  
\newcommand{\el}{\end{Lem}}  
\newcommand{\bd}{\begin{Def}\ \ }  
\newcommand{\ed}{\end{Def}}  
  
\newcommand{\pf}{\noindent{\it Proof:\ \ }}  
\newcommand{\qed}{\hfill\square}  
\newcommand{\n}{\nabla}

\newcommand{\be}{\begin{equation}}  
\newcommand{\ee}{\end{equation}}

\newcommand{\arr}{\begin{array}{rlll}}  
\newcommand{\ea}{\end{array}}  
\newcommand{\bea}{\begin{eqnarray}}  
\newcommand{\eea}{\end{eqnarray}}  
\newcommand{\bean}{\begin{eqnarray*}}  
\newcommand{\eean}{\end{eqnarray*}}  
  
\catcode`@=11  
\@addtoreset{equation}{section}  
\catcode`@=12

\begin{document}  
\begin{titlepage}  
%\rightline{hep-th/0610176}  
% \rightline{draft:  8.09.06}
~  
\vskip 1.5 true cm  
\begin{center}  
{\bigss  Flat nearly K\"ahler manifolds\\[1em]  }  
\vskip 1.0 true cm   
{\cmsslll    Vicente Cort\'es and Lars Sch\"afer} \\[3pt]  
{\tenrm   Department Mathematik, 
Universit\"at Hamburg, 
Bundesstra{\ss}e 55, 
D-20146 Hamburg, Germany\\  
cortes@math.uni-hamburg.de, schaefer@math.uni-hamburg.de}\\[1em]   
September 13, 2006
\end{center}  
\vskip 1.0 true cm  
%%%%%%%%%%%%%%%%%%%%%%%%%%%%%%%%%%%%%%%%%%%%%%%%%%%%%%%%  
\baselineskip=18pt  
\begin{abstract}  
\noindent  We classify flat strict nearly K\"ahler manifolds 
with (necessarily) indefinite metric. Any such manifold is 
locally the product of a flat  
pseudo-K\"ahler factor of maximal dimension and a strict flat nearly K\"ahler 
manifold of split signature $(2m,2m)$ with 
$m\ge 3$. Moreover, the geometry of the second factor is encoded in a
complex three-form $\zeta \in \Lambda^3 (\bC^m)^*$. 
The first nontrivial example occurs in dimension $4m=12$.\\[1em]   

\end{abstract}  
\vfill \hrule width 3.0 cm  
{\small \noindent Research of L.\ S.\ was supported by a grant of 
`Studienstiftung des deutschen Volkes'. }  
  
\end{titlepage}  
\tableofcontents  
\section{Introduction}
{\it Nearly K\"ahler geometry} originated in the study of weak holonomy groups by Gray \cite{G3}. 
In fact, nearly K\"ahler manifolds correspond 
to weak holonomy $U(n)$ and were intensively studied by Gray \cite{G0,G1,G2}. 
These manifolds appear moreover in a natural way in the Gray-Hervella classification \cite{GH} as one class of  the 
sixteen classes of almost Hermitian manifolds.

Recent interest in nearly K\"ahler manifolds came from the fact, that in dimension 6 these 
manifolds are related to the existence of Killing spinors and that they 
admit a Hermitian connection with totally 
skew-symmetric torsion. Such connections were studied by Friedrich and Ivanov 
\cite{FI} and 
are of interest in string theory. The classification 
of complete simply connected strict nearly K\"ahler manifolds 
was reduced to dimension 6 by Nagy \cite{N1,N2}.  Butruille \cite{B} 
has shown that all strict homogenous nearly K\"ahler manifolds are 3-symmetric.
These works are mainly concerned with 
Riemannian manifolds. In this paper we are especially interested in
pseudo-Riemannian metrics.

The interest in flat nearly pseudo-K\"ahler manifolds is motivated by our study \cite{S} of $tt^*$-structures (topological-antitopological fusion structures) 
on the tangent-bundle. In fact, flat nearly pseudo-K\"ahler 
manifolds provide an interesting class of $tt^*$-structures on the tangent-bundle. A second interesting class 
of solutions is given by {\it special K\"ahler manifolds} \cite{CS}. 
In other words, one can interpret $tt^*$-structures on 
the tangent-bundle as a common generalisation of these two geometries. This duality can also be seen in this work. 
We  construct  flat nearly pseudo-K\"ahler manifolds 
of split signature from a certain constant three-form, while in \cite{BC} 
flat special K\"ahler manifolds were constructed from a 
symmetric three-tensor. 

Let us describe the structure and results of the paper. 
In the first section we recall some basic facts about  
flat nearly pseudo-K\"ahler manifolds $(M,g,J)$. We give a 
self-contained presentation, with proofs which apply  
in the case of indefinite metrics and take advantage of the flatness of 
the metric. The essential points are the existence of a canonical connection 
$\n$ with skew-symmetric torsion $T$ and the properties of the tensor 
$\eta= \frac{1}{2} JDJ=D-\nabla=-\frac{1}{2}T$,  
where $D$ is the Levi-Civita connection. 

The classification results are then given in the second section. 
The first one is Theorem \ref{1stThm}, which encodes a flat 
nearly pseudo-K\"ahler structure in a constant 
three-form $\eta$ subject to two constraints. 
An explicit formula for $J$ in terms of $\eta$ is given.
Next we analyze the constraints on $\eta$. 
It turns out that the first is equivalent to require that 
$\eta$ has isotropic support (cf.\ Proposition \ref{char_i_prop}) and the 
second is equivalent to a type condition on $\eta$ (cf.\ Theorem \ref{2ndThm}).
We explicitly solve the two constraints on the real three-form 
$\eta$  (in $4m$ variables) in terms of a freely specifiable complex 
three-form $\zeta \in \Lambda^3 (\bC^m)^*$. In particular, any such form
$\zeta\neq 0$ defines a  complete simply connected 
strict flat nearly pseudo-K\"ahler manifold,  see Corollary \ref {4thCor}.  

  Further we show that any  
strict flat nearly pseudo-K\"ahler manifold is locally the 
product of a flat pseudo-K\"ahler factor 
of maximal dimension and a strict flat nearly 
pseudo-K\"ahler manifold of dimension $4m\ge 12$ and 
split signature (Theorem \ref{3rdThm}). 
This implies, in particular, 
the non-existence of strict flat nearly K\"ahler manifolds 
with positive definite metric. The work 
finishes with the classification of complete simply connected 
flat nearly K\"ahler manifolds up to isomorphism 
in terms of $GL_m(\bC)$-orbits
on $\Lambda^3 (\bC^m)^*$,  see Corollary \ref {5thCor}.

We thank Paul-Andi Nagy for useful discussions.

\section{Basic facts about flat nearly pseudo-K\"ahler manifolds}
In this section we discuss some basic properties of  nearly pseudo-K\"ahler 
manifolds. Since in this paper we are concerned with \emph{indefinite}   
nearly K\"ahler manifolds with \emph{flat} Levi-Civita connection, 
we give a self-contained discussion including indefinite metrics but 
specializing the general statements and proofs whenever possible using 
the flatness assumption. We have referred to the literature for more
general statements in the positive definite case.  

\bd
An almost complex manifold $(M,J)$ is called {\cmssl almost pseudo-Hermitian} if it 
is endowed with a pseudo-Riemannian metric $g$ which is {\cmssl pseudo-Hermitian}, i.e. 
which satifies $J^*g(\cdot,\cdot) = g(J\cdot,J \cdot)=g(\cdot,\cdot).$ 
The nondegenerate two-form $\omega := g(J \cdot , \cdot )$ is called the 
{\cmssl fundamental two-form}.\\ 
An almost pseudo-Hermitian manifold $(M,g,J)$ is called {\cmssl nearly pseudo-K\"ahler} manifold, if its Levi-Civita  connection $D$ 
satisfies the {\cmssl nearly K\"ahler condition} 
\be (D_XJ)Y = -(D_YJ)X, \quad \forall   X,Y \in \Gamma(TM). \ee
A nearly pseudo-K\"ahler manifold is called {\cmssl strict}, if $D J \ne 0.$
\ed

\bp \label{Can_con} (cf.\ Friedrich and Ivanov \cite{FI}) 
Let $(M,g,J)$ be a nearly pseudo-K\"ahler manifold. Then  there exists a unique connection $\n$ with totally skew-symmetric torsion $T^{\n}$ satisfying $\n g=0$ and $\n J =0.$  \\
More precisely, it holds 
\be T^{\nabla} = -2 \eta \mbox{ with } \eta = \frac{1}{2} J D J \ee
and $\{\eta_X,J \} =0,$ for all vector fields $X.$
\ep
To prove the proposition we give two  lemmas of independent interest.
\bl \label{lemma1}
Let $(V,g,J)$ be a pseudo-Hermitian vector space and $S$ a $(1,2)$ tensor, such that 
\begin{enumerate}
\item[(i)]
$S_X$ is skew-symmetric (with respect to $g$) for all $X \in V,$
\item[(ii)]
$T(X,Y,Z):=g(S_X  Y,Z) -g(S_Y  X,Z),$ with  $X,Y,Z \in V,$ is totally skew-symmetric.
\end{enumerate}
Then $S,$ or more precisely $(X,Y,Z) \mapsto g(S_XY,Z),$ is totally skew-symmetric, too.
\el
\pf
It holds
\bean 
g(S_XY,Z)- g(S_Y X ,Z) &\overset{(ii)}=&
-\left(g(S_Z Y,X)- g(S_Y Z ,X) \right) \\
&\overset{(i)}=&- g(S_Z Y,X)-g(S_Y X ,Z),
\eean
which implies $ g(S_X Y ,Z) = - g(S_Z Y,X). $ Together with property (i), this 
shows that $S$ is totally skew-symmetric.
\qed

\bl\label{lemma2} 
Let $(V,g,J)$ be a pseudo-Hermitian vector space and $S$ a $(1,2)$ tensor satisfying: 
\begin{enumerate}
\item[(i)]
$S$ is totally skew-symmetric and
\item[(ii)]
$[S_X,J] =0$ for all $X \in V.$
\end{enumerate}
Then $S$ vanishes.
\el
\pf With arbitrary $X,Y,Z \in V$ we show
\bean
g(S_XY,Z) &=& g(J S_XY, J Z) \overset{(ii)}= g( S_X J Y, J Z) \\
&\overset{(i)}=&-g(S_{J Z} J Y ,X)  \overset{(ii)}=-g(J S_{J Z}  Y ,X) \\
&=& g( S_{J Z}  Y , J X) \overset{(i)}= -  g( S_{Y}  J Z , J X) \\
&\overset{(ii)}=& -g(S_Y Z,X) \overset{(i)}= -g(S_XY,Z).
\eean
This shows $S=0.$
\qed \\

\pf (of Proposition \ref{Can_con})
First we show the uniqueness: \\
Let $ \n$ and $\n'$ be two such connections and $S:= \n-\n'$ their difference tensor. \\
Then from $\n J =\n'J=0$ we obtain $[S_X,J]=0$ and from   $\n g =\n'g=0$ we get the skew-symmetry of $S_X$ 
 (with respect to $g$) for all vector fields $X.$\\ 
In addition  
$$ (T^{\n} -T^{\n'})(X,Y,Z) = g(S_XY -S_YX,Z) $$
is the difference of two totally skew-symmetric tensors and hence totally skew-symmetric. 
Lemma \ref{lemma1} implies that the tensor $S$ is totally skew-symmetric and Lemma \ref{lemma2} shows the uniqueness, 
i.e. the vanishing of $S$.\\
To prove the existence we define
\bean
\n:= D- \eta \mbox{ with } \eta = \frac{1}{2} J DJ. 
\eean
The skew-symmetry of $J$ yields that $D_X J$ is skew-symmetric (with respect to $g$). Further we have
$\{ J, D_XJ\} = 0,$ as follows from deriving $J^2=-\id,$  which shows  that $\{\eta_X,J \}=0$ and that $\eta_X = JD_XJ$ is 
skew-symmetric for all vector fields $X.$ \\
{}From the skew-symmetry of $\eta_X$ and $Dg=0$ we obtain $\n g=0.$ \\
Further we compute 
\bean
\n_X J&=& D_XJ -\frac{1}{2} [JD_XJ,J] = D_XJ - \frac{1}{2}( JD_XJ J - J^2 D_XJ) \\ &=& D_XJ + J^2 D_XJ =0 \quad \forall\;  X.
\eean
This means $\n J=0.$
Finally we calculate the torsion 
\bean
T^{\n}(X,Y) &=&   D_XY -D_Y X -\eta_X Y +\eta_Y X - [X,Y] \\
&=& T^D(X,Y) - (\eta_X Y -\eta_Y X) = -\eta_X Y +\eta_Y X \\
&=& -2\eta_XY,
\eean
since $DJ$ and consequently $\eta$ is skew-symmetric by the nearly K\"ahler condition. \\
Since $\eta_X$ is skew-symmetric for all $X$ and $\eta_XY=-\eta_YX$ for all $X,Y,$ $\eta$ is totally skew-symmetric and $T^{\n}  = -2 \eta$ is totally skew-symmetric, too. 
\qed

\bp \label{eta_D_nabla_prop}
Let $(M,g,J)$ be a flat nearly pseudo-K\"ahler manifold. Then
\begin{enumerate}
\item[1)] $ \eta_X \circ \eta_Y =0$ for all $X,Y,$
\item[2)]$D \eta =\n \eta =0.$ 
\end{enumerate}
\ep

\bl (cf. Gray \cite{G0})
Let $(M,g,J)$ be a nearly pseudo-K\"ahler manifold.  Then for all $X,Y \in TM$ it is 
\bea
g(R^D(X,Y)JX,JY) - g(R^D(X,Y)X,Y)= g((D_XJ)Y,(D_XJ)Y).
\eea
\el
\pf Since in this paper we are mainly concerned with flat nearly K\"ahler manifolds, we give a short proof under the additional assumption $R^D=0.$ 
With $D$-parallel vector fields $X,Y$ we compute
\bean
0&=& g(R^D(X,Y)JX,JY) = g(D_X(D_YJ)X,JY) -  g(D_Y\underbrace{(D_XJ)X}_{=0},JY) \\
&=& Xg((D_YJ)X,JY)- g((D_YJ)X,(D_XJ)Y) \\
&=& X[ Y\underbrace{g(JX,JY)}_{const.} -g(JX,\underbrace{(D_YJ)Y}_{=0})] +g((D_XJ)Y,(D_XJ)Y)\\
&=& g((D_XJ)Y,(D_XJ)Y). \mbox{~~~~~~~~~~~~~~~~~~~~~~~~~~~~~~~~~~~~~~~~~~~~~~~~~~~~~~~~~~~~~~~~~~~~~~~~\qed}
\eean
% \qed \\
We linearize  the identity  $g((D_XJ)Y,(D_XJ)Y)=0$  in $Y$ to obtain
\bea \label{Lin_Lemma1}
g((D_XJ)Y,(D_XJ)Z) + g((D_XJ)Z,(D_XJ)Y)=0, \quad \forall X,Y,Z.
\eea
\bl (cf. Gray \cite{G1})
Let $(M,g,J)$ be a flat nearly pseudo-K\"ahler manifold. Then 
\bea
 g((D_XJ)Y,(D_ZJ)W)=0 \quad \forall X,Y,Z,W.
\eea
\el
\pf
Define  the tensor $A(X,Y,Z,W):=g((D_XJ)Y,(D_ZJ)W)=A(Z,W,X,Y).$ \\
We know $A(X,Y,X,Z)=A(X,Z,X,Y)\overset{\eqref{Lin_Lemma1}}{=}-A(X,Y,X,Z)$ for all $X,Y,Z,$ which implies
$A(X,Y,X,Z)=0$ and $A(Y,X,Z,X)=0.$ \\
We summarize the symmetries of $A:$ 
\bean
A(X,Y,Z,W)&=& -A(Y,X,Z,W)= -A(X,Y,W,Z), \text{  (nearly K\"ahler condition)}\\
A(X,Y,Z,W)&=& -A(Z,Y,X,W), \\
A(W,Y,Z,X)&=& -A(W,Y,X,Z)= A(X,Y,W,Z) =- A(X,Y,Z,W),
\eean
i.e. $A$ is totally skew-symmetric. \\
In addition it holds
\bea \label{type_of_A_I}
A(X,JY,Z,JW) &=& g((D_XJ)JY,(D_ZJ)JW) \nonumber\\&=&  g(J(D_XJ)Y,J(D_ZJ)W)=A(X,Y,Z,W), \\
A(X,Y,JZ,JW) &=& g((D_XJ)Y,(D_{JZ}J)JW) =  -  g((D_XJ)Y,(D_{W}J)J^2Z) \nonumber \\
&=&   - g((D_XJ)Y,(D_{Z}J)W)= -A(X,Y,Z,W). \label{type_of_A_II}
\eea
The skew-symmetry of $A$ yields
\bean
A(X,Y,Z,W) + A(X,Z,Y,W)    &=&0, \\
A(X,JY,Z,JW) + A(X,Z,JY,JW)&=&0.
\eean
The addition of these two equations gives
\bean 
0&=& A(X,Y,Z,W) + A(X,JY,Z,JW) + A(X,Z,Y,W) + A(X,Z,JY,JW) \\
& \overset{\eqref{type_of_A_I},\eqref{type_of_A_II}}{=}& 2 A(X,Y,Z,W) = 2g((D_XJ)Y,(D_ZJ)W)
\eean
and the lemma is proven.
\qed \\[0.8em]
\pf (of Proposition \ref{eta_D_nabla_prop})\\  1) From the last lemma we have 
\bean
0=g((D_XJ)Y,(D_ZJ)W)&=&-g((D_ZJ)(D_XJ)Y,W)\\
&=&-g(J(D_ZJ)J(D_XJ)Y,W)\\ 
&=&- 4g(\eta_Z\, \eta_X Y,W).
\eean
This shows $\eta_X \, \eta_Y =0$ for all $X,Y$ and finishes the proof of 
part 1). \\
2)  With two vector fields $X,Y$ we calculate 
\bean
(D_X\eta)_Y &=& D_X(\eta_Y)- \eta_{D_XY} \\
&\overset{D= \n +\eta}{=}&\n_X(\eta_Y) + [\eta_X ,\eta_Y] - \eta_{D_XY}\\ 
&=& (\n_X \eta)_Y + \eta_{[\n_XY -D_XY]} + [\eta_X,\eta_Y]\\
&=& (\n_X \eta)_Y - \eta_{\eta_XY } + [\eta_X,\eta_Y]\\
&\overset{n.K.}{=}& (\n_X \eta)_Y + \eta \eta_XY  + [\eta_X,\eta_Y] \overset{1)}{=} (\n_X \eta)_Y.
\eean

\noindent 
Using  $\eta= \frac{1}{2} J DJ$ and $\n J =0$ we obtain
\bean
\n \eta= \frac{1}{2} J \n( D J).  
\eean 
Therefore it is sufficient to show $  \n(DJ)=0. $
We calculate for $D$-parallel vector fields: 
\bean
g(\n_X(DJ)_YZ,W) &\overset{\n=D-\eta}{=}& g(D_X(DJ)_YZ,W)-g([\eta_X,D_YJ]Z,W)\\
& \overset{(*)}{=}&g(D_X(DJ)_YZ,W)\\
&\overset{DW=0}{=}& Xg((DJ)_YZ,W) \\
&=& X\left[ (D_Y \o)(Z,W)\right] \\
&=& D_X\left[ (D_Y \o)(Z,W)\right] \\
&\overset{DY=DZ=DW=0}=&  (D^2_{X,Y} \o)(Z,W)
,
\eean
where $\o=g(J \cdot,\cdot).$  The second term in $(*)$ vanishes by part 1), since by $\{\eta_X,J\}=0$ we get $$J[\eta_X,D_YJ]=-\{ \eta_X ,J D_YJ\} = -2\{\eta_X,\eta_Y \} =0.$$
The next lemma finishes the proof.
\qed

\bl (compare Gray \cite{G2} for the non-flat case) \\
Let $(M,g,J)$ be a flat nearly pseudo-K\"ahler manifold, then $D^2\o=0.$
\el
\pf
The nearly K\"ahler condition is equivalent to 
\be (D_X\o)(X,Y)= g((D_XJ)X,Y)=0, \; \forall X,Y. \label{Dom_equ} \ee
Further  $R^D=0$ implies the symmetry of $D_{X,Y}^2\o.$ 
Hence it suffices to show $$(D^2_{X,X} \o )(Y,Z) =0,\; \forall 
X,Y,Z.$$
Suppose $X,Y,Z$ to be $D$-parallel. Then it is
\bean
(D^2_{X,X} \o )(Y,Z)&=& D_X\left[ (D_X \o)(Y,Z) \right] \\
&\overset{\eqref{Dom_equ}}{=}&- D_X\left[ (D_Y \o)(X,Z) \right]  \\
&\overset{R^D=0}{=}& - D_Y\left[ (D_X \o)(X,Z) \right] \overset{\eqref{Dom_equ}}{=} 0.
\eean
This yields the lemma.
\qed
\section{Classification results for flat nearly pseudo-K\"ahler manifolds}
We denote by  $\bC^{k,l}$ the complex vector space $(\bC^n,J_{can})$, $n=k+l$, 
endowed with the standard $J_{can}$-invariant pseudo-Euclidian scalar
product $g_{can}$ of signature $(2k,2l)$. 

Let $(M,g,J)$ be a flat nearly pseudo-K\"ahler manifold. Then there exists for each point $p \in M$ an open set $U_p \subset M$ containing the point $p$, 
a connected open set $U_0$ of $\bC^{k,l}$ containing the origin $0\in \bC^{k,l}$ and an isometry 
$$ \Phi\,:\, (U_p,g) \tilde{\rightarrow} (U_0,g_{can}), $$
such that at the point $p$ we have: $$\Phi_* J_p = J_{can} \Phi_*.$$ 
In other words, we can suppose, that locally $M$ is a connected  open subset of $\bC^{k,l}$ containing the origin $0$ and that $g=g_{can}$ and $J_0=J_{can}.$ \\{}From Proposition \ref{Can_con} and \ref{eta_D_nabla_prop} we obtain: 
\bc \label{class_cor}  Let $M \subset  \bC^{k,l}$ be an open 
neighborhood of the origin endowed with a nearly pseudo-K\"ahler structure 
$(g,J)$ such that $g=g_{can}$ and $J_0=J_{can}$. Then  
the $(1,2)$-tensor $$ \eta := \frac{1}{2} J DJ $$ defines a constant three-form on $M \subset \bC^{k,l} =\bR^{2k,2l}$ defined by $$\eta(X,Y,Z):= g(\eta_XY,Z) $$satisfying 
\begin{enumerate}
\item[(i)] $\eta_X \, \eta_Y =0,\quad \forall X,Y,$
\item[(ii)] $\{\eta_X,J_{can}\}=0, \quad \forall X.$ 
\end{enumerate}
\ec
Conversely, we have the 
\bl \label{Lemma_J_from_eta}
Let $\eta$ be a constant three-form on an open connected 
neighborhood $M \subset \bC^{k,l}$ of $0$ satisfying (i) and (ii) of Corollary \ref{class_cor}. Then there 
exists a unique almost complex structure $J$ on $M$ such that
\begin{enumerate}
\item[a)] $J_0=J_{can},$
\item[b)] $\{\eta_X,J\}=0, \quad \forall X,$ 
\item[c)] $DJ =-2J \eta,$ 
\end{enumerate}
where $D$ stands for the Levi-Civita connection of  the 
pseudo-Euclidian vector space $\bC^{k,l}$. With 
$\n := D - \eta$ and assuming b), the last equation is equivalent to
 \begin{enumerate}
\item[c')] $\n J =0.$ 
\end{enumerate}
\el
\pf
The equivalence of c) and c') follows from a straightforward calculation. 
First we show the uniqueness of $J$:\\ 
Given two almost complex structures $J$ and $J'$ satisfying a)-c) we find 
\begin{itemize}
\item $J_0=J'_0$ and
\item $\n J =\n J' =0,$
\end{itemize} 
which shows $J=J'.$ \\
To show the existence we define
\bea
 J&=& \exp\left( 2 \sum_{i=1}^{2n} x^i\, \eta_{\partial_i}\right) J_{can} \\
&\overset{(i)}{=}& \left( Id + 2 \sum_{i=1}^{2n} x^i \, \eta_{\partial_i}\right)J_{can}, \label{JJcanEqu} 
\eea
where $x^i$ are linear coordinates of $\bC^{k,l} =\bR^{2k,2l} =\bR^{2n}$ and
$\partial_i =\frac{\partial}{\partial x^i}.$ \\
{\bf Claim:} $J$  is an almost complex structure satisfying a)-c).\\
a) From $x^i(0)= 0$ we obtain $J_0 =J_{can}.$\\
b) follows from equation \eqref{JJcanEqu} using  (i) and (ii).\\   
c) One computes 
\bean
D_{\partial_j} J &=& 2 \exp \left( 2 \sum_{i=1}^{2n} x^i \, \eta_{\partial_i} \right) \eta_{\partial_j} J_{can}\\
&\overset{(ii)}{=}& -2 \underbrace{\exp \left( 2 \sum_{i=1}^{2n} x^i \, \eta_{\partial_i} \right)\,  J_{can}}_J\, \eta_{\partial_j}=-2 J\, \eta_{\partial_j}.
\eean
It remains to prove $J^2=-Id.$ 
\bean 
J^2 &= & \left( Id + 2 \sum_{i=1}^{2n} x^i \, \eta_{\partial_i}\right)J_{can} \left( Id + 2 \sum_{i=1}^{2n} x^i \, \eta_{\partial_i}\right)J_{can}\\
&\overset{(ii)}{=}& -\left( Id + 2 \sum_{i=1}^{2n} x^i \, \eta_{\partial_i}\right) \left( Id - 2 \sum_{i=1}^{2n} x^i \, \eta_{\partial_i}\right) \\
&=& -\left[ Id -4 \left(\sum_{i=1}^{2n} x^i \, \eta_{\partial_i}\right)^2\right] \overset{(i)}{= }-Id.
\eean 
This finishes the proof.
\qed
\bt \label{1stThm} 
Let $\eta$ be a constant three-form on a connected open set $U \subset 
\bC^{k,l}$ containing $0$ which satisfies  (i) and (ii) of Corollary \ref{class_cor}. Then there exists a unique almost complex structure  
\be
J=  \exp\left( 2 \sum_{i=1}^{2n} x^i\, \eta_{\partial_i}\right) J_{can}
\ee
on $U$ such that 
\begin{enumerate}
\item[a)] $J_0 =J_{can},$ 
\item[b)] $M(U,\eta):= (U,g=g_{can},J)$ 
is a flat nearly pseudo-K\"ahler manifold.
\end{enumerate}
Any flat nearly pseudo-K\"ahler manifold is locally isomorphic to a flat 
nearly pseudo-K\"ahler manifold of the form $M(U,\eta).$ 
\et
\pf
$(M,g)$ is a flat pseudo-Riemannian manifold. Due to Lemma  
\ref{Lemma_J_from_eta}, $J$ is an almost complex structure on $M$ and $J_0=J_{can}.$ \\ 
In addition it holds
\bean
J&=& J_{can} + \left( 2 \sum_{i=1}^{2n} x^i\, \eta_{\partial_i}\right) J_{can},\eean
where $\{ \eta_{\partial_i},J \} =0$ and $ \eta_{\partial_i}$ is 
$g$-skew-symmetric. This implies that $J$ is $g$-skew-symmetric. \\
Finally from Lemma \ref{Lemma_J_from_eta}  c) and the skew-symmetry of $\eta$ it follows the skew-symmetry of $DJ.$ Therefore $(M,g,J)$ is nearly pseudo-K\"ahler. \\
The remaining statement follows from Corollary \ref{class_cor} 
and Lemma \ref{Lemma_J_from_eta}.
\qed \\

Now we discuss the general form of solutions of  
(i) and (ii) of Corollary \ref{class_cor}. In the following we shall freely 
identify the real vector space $V:=\bC^{k,l}=\bR^{2k,2l}=\bR^{2n}$ with
its dual $V^*$ by means of the pseudo-Euclidian scalar product $g=g_{can}$. 

\bp \label{char_i_prop}
A three-form $\eta \in \Lambda^3 V^*\cong \Lambda^3 V$ 
satisfies (i) of Corollary \ref{class_cor} if and only if 
there exists  an isotropic subspace $L\subset V$
such that $\eta \in \Lambda^3 L \subset \Lambda^3 V$. 
If $\eta$ satisfies (i) and (ii) of Corollary \ref{class_cor} then
there exists a $J_{can}$-invariant isotropic subspace $L\subset V$ 
with $\eta \in \Lambda^3 L$. 
\ep

\noindent 
{\bf Remark:} From the Proposition \ref{char_i_prop} we conclude that there 
are no strict flat nearly pseudo-K\"ahler manifolds of  dimension 
less than 8. We shall see later that the  dimension cannot be smaller 
than 12, see Corollary \ref {Dim6Cor}. 

\noindent 
We define the {\cmssl support} of $\eta \in \Lambda^3V$ by
\be \Sigma_{\eta}:= {\rm span}\{\eta_XY \,|\, X,Y \in V\} \subset V.\ee

\pf (of Proposition \ref{char_i_prop}) The proposition follows from the next two lemmas by taking $L=\Sigma_\eta$. 
\qed 

\bl \label{Lemma_cond_I}
$\Sigma_{\eta}$ is isotropic if and only if $\eta$ satisfies (i) of  
Corollary \ref{class_cor}. If $\eta$ satisfies (ii) of  Corollary 
\ref{class_cor}, then
$\Sigma_{\eta}$ is $J_{can}$-invariant. 
\el

\pf
First the isotropy of $\Sigma_{\eta}$ is equivalent to $g(\eta_XY,\eta_ZW)=0$ for all $X,Y,Z,W \in V.$  \\
Further it holds 
\be g(\eta_XY,\eta_ZW) \overset{(*)}{=} -g(\eta_Z \eta_XY,W). \label{calc} \ee
In $(*)$ we used 
$$ g(\eta_XY,Z)= \eta(X,Y,Z) =-  \eta(X,Z,Y) =-g(\eta_X Z, Y) = - g(Y,\eta_X Z) , \; \forall X,Y,Z. $$
Equation \eqref{calc}  shows  $\eta_X \, \eta_Y =0$ for all $X,Y \in V$ if and only if $\Sigma_{\eta}$ is isotropic. The last assertion follows from
\[ J\Sigma_{\eta} = {\rm span}\{J\eta_XY \,|\, X,Y \in V\} \overset{(ii)}{=} 
{\rm span}\{-\eta_XJY \,|\, X,Y \in V\} = \Sigma_{\eta}.\] 
\qed \\

\bl Let $\eta \in \Lambda^3V$. 
Then $\eta \in \Lambda^3\Sigma_{\eta}.$ 
\label{Lemma_val_insubspace}  \el 
\pf
We take a complement  $W\subset V$ of $L=\Sigma_{\eta}.$  The decomposition 
$$ \Lambda^3V = \underset{p+q=3}{\bigoplus} \Lambda^p L \wedge \Lambda^q W$$
induces a decomposition 
$$ \eta= \sum_{p+q=3} \eta^{p,q}.$$
Taking $X,Y \in   L^{\bot}$ yields $L \ni \eta_XY = \eta^{0,3}_XY + \eta^{1,2}_XY.$ Now since $\eta^{0,3}_XY \in W$ and  $\eta^{1,2}_XY \in L$, we get 
$\eta^{0,3}=0.$ Further the choice $X \in L^{\bot}$ and $Y \in W^{\bot}$ yields $\eta^{1,2}=0$ and then the choice $X,Y \in   W^{\bot}$ yields $\eta^{2,1}=0.$ 
This shows $\eta = \eta^{3,0}.$
\qed \\

Any three-form $\eta$ on $(V,J_{can})$ decomposes with respect to the 
grading induced by the decomposition
$$ V_{\bC} = V^{1,0} \oplus V^{0,1} $$
into 
\be \eta = \eta^+ + \eta^- \ee 
with 
$$\eta^+ \in \Lambda^+V := (\Lambda^{2,1}V + \Lambda^{1,2}V)^{\rho}$$
and 
$$\eta^- \in  \Lambda^-V := (\Lambda^{3,0}V + \Lambda^{0,3}V)^{\rho},$$
where $\rho$ is the canonical real structure on $V_{\bC}$ with real-points $V$ 
which extends to the exterior algebra.
\bt \label{2ndThm} 
A three-form $\eta \in \Lambda^3 V^*\cong \Lambda^3 V$ satisfies (i)
and (ii) of Corollary \ref{class_cor} if and only if there exists  an 
isotropic $J_{can}$-invariant subspace $L\subset V$ such that 
$\eta \in \Lambda^-L=(\Lambda^{3,0}L + \Lambda^{0,3}L)^{\rho}
\subset \Lambda^3L\subset \Lambda^3 V$. (The smallest such subspace $L$ is 
$\Sigma_\eta$.)  
\et 

\pf By Proposition \ref{char_i_prop}, the conditions (i)
and (ii) of Corollary \ref{class_cor} imply the existence
of an isotropic $J_{can}$-invariant subspace $L\subset V$ such that 
$\eta \in \Lambda^3L$. The next lemma shows that the condition (ii) 
is equivalent to $\eta \in \Lambda^-V.$ Therefore 
$\eta \in \Lambda^3L\cap  \Lambda^-V = \Lambda^-L$. The converse statement
follows from the same argument. \qed 

\bl 
It is 
$$  (\Lambda^{3,0}V + \Lambda^{0,3}V)^{\rho}= \{ \eta \in \Lambda^3 V\,|\, 
\eta(\cdot,J\cdot,J\cdot) = -\eta(\cdot,\cdot,\cdot) \} = \{ \eta \in \Lambda^3 V\,|\,\{\eta_X,J\}=0, \; \forall X \in V \}.$$
\el
\pf We have the decomposition 
\[ \Lambda^2V = (\Lambda^{1,1}V)^\rho \oplus (\Lambda^{2,0}V + 
\Lambda^{0,2}V)^\rho,\]
where 
\[(\Lambda^{1,1}V)^\rho  = \{ \a \in \Lambda^2V| \a (J \cdot , J \cdot )= \a\} 
\cong \{ A\in \mathfrak{so}(V)| [A,J] = 0\}\] 
and 
\[(\Lambda^{2,0}V +\Lambda^{0,2}V )^\rho  = \{ \a \in \Lambda^2V| 
\a (J \cdot , J \cdot )= -\a\} 
\cong \{ A\in \mathfrak{so}(V)| \{A,J\} = 0\}.\]
This induces the following direct decomposition:  
\[ V\otimes \Lambda^2V = V\otimes(\Lambda^{1,1}V)^\rho + V\otimes
(\Lambda^{2,0}V + 
\Lambda^{0,2}V)^\rho.\] 
We claim that  
\[ (V\otimes
 (\Lambda^{2,0}V + 
\Lambda^{0,2}V)^\rho) \cap \Lambda^3V = (\Lambda^{3,0}V
+ \Lambda^{0,3}V)^\rho.\]
The claim implies the lemma. To see the claim, let us first observe the
following obvious inclusion:
\[ (V\otimes
 (\Lambda^{2,0}V + 
\Lambda^{0,2}V)^\rho) \cap \Lambda^3V \supset (\Lambda^{3,0}V
+ \Lambda^{0,3}V)^\rho.\] 
To show the equality 
we observe that an element of $V\otimes
 (\Lambda^{2,0}V + \Lambda^{0,2}V)^\rho$ is totally skew if and only if
its four components in 
\[ V^{1,0}\otimes \Lambda^{2,0}V,\quad  V^{1,0}\otimes \Lambda^{0,2}V,\quad
V^{0,1}\otimes \Lambda^{2,0}V\quad \mbox{and}\quad V^{0,1}\otimes  
\Lambda^{0,2}V\]
are totally skew.  To finish we notice that
\[ (V^{1,0}\otimes \Lambda^{2,0}V +  V^{0,1}\otimes  
\Lambda^{0,2}) \cap \Lambda^3 V = (\Lambda^{3,0}V + \Lambda^{0,3}V)^{\rho}\]
and  
\[ (V^{1,0}\otimes \Lambda^{0,2}V + 
V^{0,1}\otimes \Lambda^{2,0})\cap \Lambda^3 V = 0.\]  
\qed

\bc \label{Dim6Cor} There are no strict  flat nearly pseudo-K\"ahler 
manifolds of 
dimension 
less than 12.
\ec 

\pf By Theorem \ref{1stThm} and \ref{2ndThm} any  flat nearly pseudo-K\"ahler 
manifold $M$ 
is locally of the form $M(U,\eta )$, where $\eta \in \Lambda^-L$ for 
an isotropic $J_{can}$-invariant subspace $L\subset V$ and $U\subset V$ is 
an open subset. $M(U,\eta )$ 
is strict if and only if $\eta\neq 0$, which is possible only for
$\dim_\bC L \ge 3$, i.e.\ for $\dim M\ge 12$. \qed

\bt \label{3rdThm} Any strict flat nearly pseudo-K\"ahler 
manifold is locally a pseudo-Riemannian product $M=M_0\times M(U,\eta )$ of a 
flat pseudo-K\"ahler factor $M_0$ of maximal dimension and a 
strict flat nearly pseudo-K\"ahler  
manifold $M(U,\eta )$ of (real) signature $(2m,2m)$, 
$4m=\dim M(U,\eta )\ge 12$. The  $J_{can}$-invariant isotropic support 
$\Sigma_\eta$ has complex dimension  $m$. 
\et 

\bc
Let $(M,g,J)$ be a flat nearly K\"ahler manifold with a  
(positive or negative) definite metric $g$ then 
$\eta =0,$ 
$\n=D$  and $DJ=0,$ i.e. $(M,g,J)$ is a K\"ahler manifold. 
\ec

\pf (of Theorem \ref{3rdThm})  By Theorem \ref{1stThm} and \ref{2ndThm}, 
$M$ is locally isomorphic to an open subset of a manifold of 
the form $M(V,\eta )$, where $\eta \in \Lambda^3V$ has a $J_{can}$-invariant 
and isotropic support $L=\Sigma_\eta$. We choose a   $J_{can}$-invariant
isotropic subspace $L'\subset V$ such that $V':=L+L'$ is nondegenerate and 
$L\cap L'=0$ and put $V_0 = (L+L')^\perp$. Then $\eta \in \Lambda^3V'
\subset  \Lambda^3V$  and $M(V,\eta ) = M(V_0,0)\times M(V',\eta )$. 
Notice that $M(V_0,0)$ is simply the flat pseudo-K\"ahler manifold $V_0$ 
and that  $M(V',\eta )$ is strict and of split signature $(2m,2m)$, where 
$m=\dim_\bC L\ge 3$. \qed 
   
For the rest of this paper we consider the case $V\cong\bC^{m,m}$ 
and denote a maximal $J_{can}$-invariant isotropic 
subspace by $L$. We will say that a  complex three-form
$\zeta \in \Lambda^3(\bC^m)^*$ has \emph{maximal support} 
if ${\rm span} \{ \zeta (Z,W,\cdot )| Z,W\in \bC^m\}
=(\bC^m)^*$. 

\bc \label{4thCor} 
 Any non-zero complex three-form $\zeta \in \Lambda^{3,0}L\cong 
\Lambda^3(\bC^m)^*$ defines
a complete flat  simply connected strict nearly pseudo-K\"ahler  manifold 
$M(\eta):= M(V,\eta )$, $\eta = \zeta + \bar\zeta \in \Lambda^3L \subset 
\Lambda^3V$, of split signature. 
$M(\eta)$ has no pseudo-K\"ahler de Rham factor 
if and only if $\zeta$ has maximal support. 

Conversely,
any complete flat  simply connected nearly pseudo-K\"ahler  manifold
without pseudo-K\"ahler de Rham factor is of this form. 
\ec 

\pf This follows from the previous results observing that the 
support of $\eta$ is maximally isotropic if and only if $\zeta$ has 
maximal support. \qed  

\bc \label{5thCor} 
The map $\zeta \mapsto M(\zeta + \bar\zeta)$ induces  a bijective 
correspondence between $GL_m(\bC)$-orbits on 
the open subset $\Lambda^3_{reg}(\bC^m)^*\subset \Lambda^3(\bC^m)^*$ 
of three-forms $\zeta$ with maximal support and isomorphism classes of 
complete flat  simply connected nearly pseudo-K\"ahler  manifolds 
$M(\zeta + \bar\zeta)$ of real dimension $4m\ge 12$ and 
without pseudo-K\"ahler de Rham factor. 
\ec

\end{document}